\documentclass{amsart}
\usepackage{amssymb}
\usepackage{amsfonts}
\usepackage{pstricks}

\theoremstyle{plain}

\renewcommand\bigskip\medskip

\def\F{\mathbb{F}}

\def\Q{\mathbb{Q}}
\def\leq{\leqslant}
\def\geq{\geqslant}

\begin{document}
\centerline{\bf{ON A TWO-VALUED INFINITE SEQUENCE}}
\centerline{\bf{ AND A RELATED CONTINUED FRACTION IN $\Q((T^{-1}))$}}

\vskip 0.5 cm
\centerline{\bf{by A. Lasjaunias }}
\centerline {(Bordeaux, France)}
\vskip 0.5 cm

{\bf{Abstract.}} This note is a complement to an article which was published, six years ago, in The Ramanujan Journal (vol. 45.3, 2018). Here, the goal is to fully describe a singular transcendental continued fraction in $\Q((T^{-1}))$, tied to a particular infinite two letters word.

\vskip 0.5 cm
\par In a previous paper [1], a particular infinite word over the alphabet $\{1,2\}$ was considered, leading to a generating function in $\Q((T^{-1}))$. The continued fraction expansion of this function was considered and it could be described partially.  At the end of the article [1, p.870--871],  a full description of this continued fraction was stated as a conjecture. The aim of this note is to complete the article [1], by proving the conjectured formulas stated there. 
\par We first recall the definition of this infinite word. Let $(W_{n})_{n\geqslant 0}$ be the sequence of finite words over the
alphabet $\{{1,2\}}$, defined recursively as follows:
$$W_{0}=\emptyset ,\text{ }W_{1}=1,\quad \text{ and }\quad 
W_{n}=W_{n-1},2,W_{n-2},2,W_{n-1},\quad \text{ for }\quad n\geqslant 2.$$
Then the infinite word $W$ will be the projective limit of the sequence $(W_{n})_{n\geqslant 0}$, that is the word beginning by $W_n$ for all $n\geq 0$. Hence, we have
$$W=(w_i)_{i\geq 1}=122121212212....$$
We must indicate that this infinite word is derived from an example of a continued fraction expansion for a quartic power series over $\F_3$, introduced by Mills and Robbins in 1986 ([2, p. 403]). For further explanations and considerations about this continued fraction and this infinite word, the reader may consult [3, p. 226--227]) and [4, p. 57--58]. Eventhough this infinite sequence is very simply defined, it was proved in [4] that it is not an automatic sequence. 

\par   We associate to this infinite word $W=(w_i)_{i\geq 1}$, a generating function $\theta \in \Q((T^{-1}))$, defined by  
$$\theta=\sum_{i=1}^{\infty} w_iT^{-i}.$$
We are concerned with the continued fraction expansion of $\theta$ in the field $\Q((T^{-1}))$. For a basic introduction to continued fractions, particularly in power series fields, the reader may consult [5], and also [3] for a deeper and more general account concerning power series fields. Here, we prove the following theorem.  

\vskip 0.5 cm
\noindent {\bf{Theorem}} \emph{ Let $W$ and $\theta$ be defined as above. Then $\theta$ is an irrational element of $\Q((T^{-1}))$ which can be expanded as an infinite continued fraction. We have : $$\theta=[0,a_1,a_2,\dots,a_n,\dots]\in \Q((T^{-1}))$$ where the partial quotients, $a_n\in \Q[T]$, are non constant polynomials. The first four elements of this sequence of partial quotients $(a_n)_{n\geq 1}$ are given by 
  $$\eqno{(E_0)}\quad (a_1,a_2,a_3,a_4)=(T-2, T/2+1/4, 8T/5+76/25, -125T/48+25/24).$$
 To describe the sequence $(a_n)_{n\geq 5}$, we need to introduce the following elements. The first one is a sequence of positive integers, $(\ell_n)_{n\geq 0}$, defined by $$\ell_0=0,\quad \ell_1=1\quad \textrm{ and}\quad \ell_{n+1}=2\ell_n+\ell_{n-1}+2 \quad \textrm{ for}\quad n\geq 1.$$
 (Note that $\ell_n$ is simply the length of the word $W_n$ for $n\geq 0$). 
 Then we introduce two sequences of rational numbers, $(r_n)_{n\geq 1}$ and $(s_n)_{n\geq 1}$. For $n\geq 1$, we set :
 $$r_n=(4/25)(2\ell_n-\ell_{n-1}+1)\quad \text{and} \quad s_n=r_{n+1}+r_n.$$
 Based upon these last sequences, we consider four more sequences of rational numbers $(\lambda_{4n+i})_{n\geq 1}$, for $i=1,2,3$ and $4$, which are defined as follows, for $n\geq 1$:
$$
\begin{array}{ll}
\lambda_{4n+1}=(-1)^{n+1}r_n^2, & \qquad \lambda_{4n+2}=(-1)^{n+1}(r_ns_n)^{-1},\\
\lambda_{4n+3}=(-1)^{n+1}s_n^2,& \qquad \lambda_{4n+4}=(-1)^{n+1}(r_{n+1}s_n)^{-1}.
\end{array} \eqno{(L_n)}
$$
Finally, for $n\geq 1$,  we introduce two polynomials, $A_n$ and $B_n$ in $\Q[T]$, defined by
 $$A_n=(T^{(3\ell_n+\ell_{n-1}+3)/2}+T^{(\ell_n+\ell_{n-1}+1)/2}-2)/(T-1)$$
  and $$B_n=(T^{(\ell_n+\ell_{n-1}+3)/2}-1)/(T-1).$$
 The description of the sequence $(a_n)_{n\geq 5}$ is now completed by giving the following formulas $(E_n)$, for $n\geq 1$ :
 $$(a_{4n+1},a_{4n+2},a_{4n+3},a_{4n+4})=(\lambda_{4n+1}A_n,\lambda_{4n+2}(T-1),\lambda_{4n+3}B_n,\lambda_{4n+4}(T-1)).$$}
\par {\bf{Remark :}} As a consequence of this theorem, we can compute the value of the irrationality measure for $\theta$. Indeed, for any irrational element in a power series field, represented by an infinite continued fraction, the value of its irrationality measure is directly connected to the sequence of the degrees of the partial quotients. For an introduction on this matter, the reader may consult [5, p. 11-13]. In the present case, the knowledge of the infinite continued fraction for $\theta$ implies that this irrationality measure is equal to 3 (see [1, p. 863]). It is known that Roth Theorem, on rational approximation of algebraic real numbers, has an analogue in power series fields over a field of characteristic zero. Hence, all  elements in $\Q((T^{-1}))$, algebraic over $\Q(T)$, have an irrationality measure equal to 2. Therefore  $\theta \in \Q((T^{-1}))$ is a transcendental element.  
\vskip 1 cm
\par {\bf{Proof of the theorem:}} First, we observe that the irrationality of $\theta$ in $\Q((T^{-1}))$ has been proved in [1, Theorem 3, p. 862] (this follows from arguments recalled below). We use the definitions and notations, concerning continued fractions, convergents and continuants, as presented in [5, p. 1-7]. Hence, by truncation of the infinite continued fraction for $\theta$, we have :
$$x_n/y_n=[0,a_1,a_2,...,a_n]=1/[a_1,...,a_n], \quad \text{for} \quad n\geq 1. \eqno{(1)}$$
For $n\geq 1$, $x_n/y_n$ is a convergent of theta in $\Q(T)$, and $x_n$ and $y_n$ are coprime polynomials in $\Q[T]$, called continuants, built from the partial quotients $a_i\in \Q[T]$, and expressed by
$$x_n=\langle a_2,a_3,...,a_n \rangle \quad \text{and} \quad y_n=\langle a_1,a_2,...,a_n \rangle ,\quad \text{for} \quad n\geq 1. \eqno{(2)}$$  
(Note that the empty continuant is equal to $1$. Hence $x_1/y_1=1/a_1$, $x_2/y_2=a_2/(a_1a_2+1)$,  $x_3/y_3=(a_2a_3+1)/(a_1a_2a_3+a_1+a_3)$,  etc... ).
Indeed, in the present case ($a_0=0$), these continuants are both defined by the same following recursive relation: 
$$K_n=a_nK_{n-1}+K_{n-2} \quad \text{with} \quad (x_1,x_0)=(1,0) \quad \text{and} \quad (y_1,y_0)=(a_1,1). \eqno{(3)}$$  
For $1\leq m <n$, we set $\Delta(m,n)=x_ny_m-x_my_n$. Adapting the formula stated in [5, p. 7], for $1\leq m <n$, we have the following
$$\Delta(m,n)=x_ny_m-x_my_n=(-1)^m \langle a_{m+2},a_{m+3},...,a_n \rangle.\eqno{(4)}$$
(Note that for $m=n-1$, the continuant on the right side of $(4)$ is empty, and this becomes simply the classical formula: $x_ny_{n-1}-x_{n-1}y_n=(-1)^{n-1}$.)
\par Now, we shall report several results, concerning the infinite continued fraction for $\theta$, which have been obtained during the proof of Theorem 3 in [1]. There are two sequences of rational elements $(R_n/S_n)_{n\geq 1}$ and $(R'_n/S'_n)_{n\geq 1}$
which are very good rational approximations to $\theta$, and therefore convergents of $\theta$. (The existence of these rational approximations is a direct proof of the irrationality of $\theta$). The four polynomials $R_n$, $S_n$, $R'_n$ and $S'_n$ are unitary and pairwise coprime. The polynomials $S_n$ and $S'_n$ are simply described using the sequence $(\ell_n)_{n\geq 0}$. Furthermore, for $n\geq 1$, the rank of these convergents is $4n$ for $R_n/S_n$ and $4n+2$ for $R'_n/S'_n$. All this follows from Lemma 1, Lemma 2, Lemma 3 and the last part of the proof of Theorem 3 in [1]. To report all this, we need to use a particular notation.\newline For $n\geq 1$, the leading coefficient of the partial quotient $a_n$ is denoted by $\lambda_n$. The leading coefficient of $x_n$ and $y_n$, following from $(2)$ and $(3)$, is obtained as the product of the leading coefficients of the partial quotients involved in its first term. For $n\geq 1$, we set $\mu_n=\prod_{i=1}^{n}\lambda_i$. 
Since we have $a_1=T-2$ (as will be shown below), we see that both $x_n$ and $y_n$ have $\mu_n$ as leading coefficient. For $n\geq 1$, we introduce a pair of unitary polynomials $(x^*_n,y^*_n)$ linked to the pair $(x_n,y_n)$, by defining 
$$x_n=\mu_n x^*_n \quad \text{and} \quad y_n=\mu_n y^*_n  \quad \text{for} \quad n\geq 1. \quad \eqno{(5)}$$ 
With this notation, the results obtained in the previous article [1], and mentioned above, can be resumed by the following four statements ($(6)$ to $(9)$):
$$x^*_{4n}=R_n \quad \text{and}\quad  y^*_{4n}=S_n=T^{(\ell_n+\ell_{n-1}+3)/2}(T^{\ell_n+1}-1) \quad \text{for} \quad n\geq 1. \eqno{(6)}$$
$$x^*_{4n+2}=R'_n \quad \text{and}\quad  y^*_{4n+2}=S'_n=T^{3\ell_n+\ell_{n-1}+4}-1 \quad \text{for} \quad n\geq 1. \eqno{(7)}$$
and
$$R_nS'_n-R'_nS_n=S_{n+1}R'_n-R_{n+1}S'_n=(-1)^n(T-1)\quad \text{for} \quad n\geq 1. \eqno{(8)}$$
Here above, in $(6)$ and $(7)$, we have given the simple form of the polynomials $S_n$ and $S'_n$, while the form of $R_n$ and $R'_n$ is more complicated. A construction of these last polynomials is explained in [1]. However, here, we will only need the following [1, p. 868]
$$R_1=T^3+2T^2+T-1\quad \text{and}\quad R'_1=T^6+2T^5+2T^4+T^3+2T^2+T+2. \eqno{(9)}$$
\par First, we shall obtain the form of the partial quotients $a_{4n+2}$ and $a_{4n+4}$ for $n\geq 1$, by combining $(6)$, $(7)$, $(8)$ and the formula $(4)$. Since $\langle x\rangle =x$ in $\Q[T]$, $(4)$ implies
$$\Delta(4n,4n+2)=a_{4n+2} \quad \text{and}\quad \Delta(4n+2,4n+4)=a_{4n+4}. \eqno{(10)}$$ 
Using $(4)$, $(5)$, $(6)$, $(7)$ and $(8)$, for $n\geq 1$, we can write
$$\Delta(4n,4n+2)=\mu_{4n+2}R'_n\mu_{4n}S_n-\mu_{4n+2}S'_n\mu_{4n}R_n=\mu_{4n+2}\mu_{4n}(-1)^{n+1}(T-1). \eqno{(11)}$$ 
In the same way, for $n\geq 1$, we also have
$$\Delta(4n+2,4n+4)=\mu_{4n+4}R_{n+1}\mu_{4n+2}S'_n-\mu_{4n+4}S_{n+1}\mu_{4n+2}R'_n$$ 
$$\Delta(4n+2,4n+4)=\mu_{4n+2}\mu_{4n+4}(-1)^{n+1}(T-1). \eqno{(12)}$$
Hence, combining $(10)$, $(11)$ and $(12)$, for $n\geq 1$, we have
 $$a_{4n+2}=(-1)^{n+1}\mu_{4n+2}\mu_{4n}(T-1),\qquad  a_{4n+4}=(-1)^{n+1}\mu_{4n+2}\mu_{4n+4}(T-1). \eqno{(13)}$$ 
 From formulas $(13)$, for $n\geq 1$, we obtain 
 $$\lambda_{4n+2}=\mu_{4n+2}\mu_{4n}(-1)^{n+1}\quad \text{and}\quad \lambda_{4n+4}=\mu_{4n+2}\mu_{4n+4}(-1)^{n+1}.\eqno{(14)}$$
  \par For $n\geq 1$, let us now compute $a_{4n+3}$. To do so , we will need first the following equality :
 $$S_n+S_{n+1}=S'_n[(T-1)B_n+1]\quad \text{for} \quad n\geq 1.\eqno{(15)}$$
 This equality shows the link between the polynomials $S_n$,$S_{n+1}$, $S'_n$, introduced in $(6)$ and $(7)$, and the polynomial $B_n$ introduced in the theorem. The rightness of $(15)$ is obtained directly by a basic computation.
 Now, to obtain $a_{4n+3}$, we will use $\Delta(4n,4n+4)$. Recalling that $\langle x,y,z\rangle=xyz+x+z$ in $\Q[T]$, appplying $(4)$, for $n\geq 1$, we have 
 $$\Delta(4n,4n+4)=a_{4n+2}a_{4n+3}a_{4n+4}+a_{4n+2}+a_{4n+4}.$$
 With $(13)$ and $(14)$, for $n\geq 1$, this becomes
 $$\Delta(4n,4n+4)=\lambda_{4n+2}\lambda_{4n+4}a_{4n+3}(T-1)^2+(\lambda_{4n+2}+\lambda_{4n+4})(T-1). \eqno{(16)}$$ 
 On the other hand, using $(4)$ $(5)$, $(6)$ and $(7)$, for $n\geq 1$, we can write
 $$\Delta(4n,4n+4)=\mu_{4n}\mu_{4n+4}(R_{n+1}S_n-S_{n+1}R_n). \eqno{(17)}$$ 
 We need to transform $(17)$. Indeed, applying $(8)$, we can write 
  $$R_{n+1}/S_{n+1}-R_n/S_n=R_{n+1}/S_{n+1}-R'_n/S'_n+R'_n/S'_n-R_n/S_n$$
  $$=(-1)^{n+1}(T-1)(1/S_{n+1}+1/S_n)/S'_n. \eqno{(18)}$$ 
 Using $(15)$, $(18)$ implies
 $$R_{n+1}S_{n}-R_nS_{n+1}=(-1)^{n+1}(T-1)(S_{n}+S_{n+1})/S'_n$$
   $$=(-1)^{n+1}(T-1)[(T-1)B_n+1]. \eqno{(19)}$$ 
 In order to get $a_{4n+3}$, we have to compare formulas $(16)$ and $(17)$.  Applying $(19)$, for $n\geq 1$, we get  $$\lambda_{4n+2}\lambda_{4n+4}a_{4n+3}(T-1)+(\lambda_{4n+2}+\lambda_{4n+4})=$$
 $$(-1)^{n+1}\mu_{4n}\mu_{4n+4}[(T-1)B_n+1]. \eqno{(20)}$$
 From $(14)$, for $n\geq 1$, we observe the following 
 $$(\lambda_{4n+2}\lambda_{4n+4})/(\mu_{4n}\mu_{4n+4})=\mu_{4n+2}^2 \eqno{(21)}$$
 and 
 $$(\lambda_{4n+2}+\lambda_{4n+4})/(\mu_{4n}\mu_{4n+4})=(-1)^{n+1}\mu_{4n+2}(1/\mu_{4n}+1/\mu_{4n+4}).\eqno{(22)}$$
 Using $(21)$ and $(22)$, $(20)$ becomes
 $$\mu_{4n+2}^2(T-1)a_{4n+3}=(-1)^{n+1}[(T-1)B_n+1-\mu_{4n+2}(1/\mu_{4n}+1/\mu_{4n+4})].\eqno{(23)}$$
 Here, we come to a critical point. Both sides of $(23)$ must vanish for $T=1$ and therefore we must have
 $$1/\mu_{4n+2}=1/\mu_{4n}+1/\mu_{4n+4}\quad \text{for} \quad n\geq 1.\eqno{(III)}$$
 Consequently, $(23)$ reduces to 
 $$\mu_{4n+2}^2a_{4n+3}=(-1)^{n+1}B_n.\eqno{(24)}$$
  Hence, we have
 $$a_{4n+3}=\lambda_{4n+3}B_n\quad \text{with}\quad \lambda_{4n+3}=(-1)^{n+1}\mu_{4n+2}^{-2}\quad \text{for}\quad n\geq 1.\eqno{(25)}$$
 
  \par Finally, we turn to the partial quotient $a_{4n+1}$. We will follow the same lines as in the previous case. For technical reasons, we have to consider $n\geq 2$. In this case, we will need another particular equality which is the following :
  $$S'_n-S'_{n-1}=S_n[(T-1)A_n+2]\quad \text{for} \quad n\geq 2.\eqno{(26)}$$
  This equality shows the link between the polynomials $S'_n$,$S'_{n-1}$, $S_n$, introduced in $(6)$ and $(7)$, and the polynomial $A_n$ introduced in the theorem. Again, the proof of $(26)$ is obtained directly by a basic computation.
  Now, in order to obtain $a_{4n+1}$, we will use $\Delta(4n-2,4n+2)$. Recalling that $\langle x,y,z\rangle=xyz+x+z$, appplying $(4)$, for $n\geq 1$, we have 
  $$\Delta(4n-2,4n+2)=a_{4n}a_{4n+1}a_{4n+2}+a_{4n+1}+a_{4n+2}.$$
  Introducing $(13)$ and $(14)$, for $n\geq 2$, this becomes
  $$\Delta(4n-2,4n+2)=\lambda_{4n}\lambda_{4n+2}a_{4n+3}(T-1)^2+(\lambda_{4n}+\lambda_{4n+2})(T-1). \eqno{(27)}$$ 
  On the other hand, using $(4)$, $(5)$, $(6)$ and $(7)$, for $n\geq 2$, we can write
  $$\Delta(4n-2,4n+2)=\mu_{4n-2}\mu_{4n+2}(R'_{n}S'_{n-1}-S'_{n}R'_{n-1}). \eqno{(28)}$$ 
  We need to transform $(28)$. As above, using $(8)$, for $n\geq 2$, we can write
  $$R'_{n}/S'_{n}-R'_{n-1}/S'_{n-1}=R'_{n}/S'_{n}-R_n/S_n+R_n/S_n-R'_{n-1}/S'_{n-1}$$
   $$=(-1)^{n+1}(T-1)(1/S'_{n-1}-1/S'_n)/S_n. \eqno{(29)}$$ 
   Applying $(26)$, $(29)$ implies
  $$R'_{n}S'_{n-1}-S'_{n}R'_{n-1}=(-1)^{n}(T-1)(S'_{n}-S'_{n-1})/S_n$$
    $$=(-1)^{n}(T-1)[(T-1)A_n+2]. \eqno{(30)}$$ 
    Now, we compare $(27)$ and $(28)$. Using $(30)$, we get, for $n\geq 2$
    $$\lambda_{4n}\lambda_{4n+2}a_{4n+1}(T-1)+(\lambda_{4n}+\lambda_{4n+2})=(-1)^{n}\mu_{4n-2}\mu_{4n+2}[(T-1)A_n+2]. \eqno{(31)}$$
    From $(14)$, for $n\geq 2$, we observe the following 
     $$(\lambda_{4n}\lambda_{4n+2})/(\mu_{4n-2}\mu_{4n+2})=-\mu_{4n}^2 \eqno{(32)}$$
     and 
     $$(\lambda_{4n}+\lambda_{4n+2})/(\mu_{4n-2}\mu_{4n+2})=(-1)^{n}\mu_{4n}(1/\mu_{4n+2}-1/\mu_{4n-2}).\eqno{(33)}$$   
     Using $(32)$ and $(33)$, $(31)$ becomes
     $$\mu_{4n}^2(T-1)a_{4n+1}=(-1)^{n+1}[(T-1)A_n+2-\mu_{4n}(1/\mu_{4n+2}-1/\mu_{4n-2})].\eqno{(34)}$$
     As above, this last equality has an important consequence. Both side of $(35)$ must vanish for $T=1$ and therefore we must have
     $$1/\mu_{4n+2}-1/\mu_{4n-2}=2/\mu_{4n}\quad \text{for} \quad n\geq 2.\eqno{(I)}$$
     Consequently, $(34)$ reduces to 
        $$\mu_{4n}^2a_{4n+1}=(-1)^{n+1}A_n.\eqno{(35)}$$
         Hence, we have
        $$a_{4n+1}=\lambda_{4n+1}A_n\quad \text{with}\quad \lambda_{4n+1}=(-1)^{n+1}\mu_{4n}^{-2}\quad \text{for}\quad n\geq 2.\eqno{(36)}$$ 
   \par We will now look at the first partial quotients, using the pairs $(R_1,S_1)$ and  $(R'_1,S'_1)$. According to $(1)$, $(6)$ and $(7)$, we have 
   $$R_1/S_1=x_4/y_4=[0,a_1,a_2,a_3,a_4]\quad \text{and}\quad R'_1/S'_1=x_6/y_6=[0,a_1,a_2,a_3,a_4,a_5,a_6].$$ We take into account the values for $(R_1,R'_1)$, given in $(9)$, and also, from $(7)$ and $(8)$, $(S_1,S'_1)=(T^4-T^2,T^7-1)$. By expanding the rational functions $R_1/S_1$ and $R'_1/S'_1$, we get the first four partial quotients, as they are declared by $(E_0)$ in the theorem, but also 
   $$ (a_5,a_6)= ((144/625)(T^2 +T + 2), (625/528)(T -1)).$$      
   This shows that $(36)$ holds also for $n=1$. Moreover, we have $\mu_4=-25/12$ and, applying $(14)$ and $(III)$ for $n=1$, from $\mu_6$, we get $\mu_8=-25/32$. We introduce two sequences of rational numbers  $(r_n)_{n\geq 1}$ and  $(s_n)_{n\geq 1}$,
   defined by
   $$r_n=-\mu_{4n}^{-1}\quad \text{and}\quad s_n=-\mu_{4n+2}^{-1} \quad \text{for}\quad n\geq 1.$$
   According to $(14)$, $(25)$ and $(36)$, for $n\geq 1$, $(E_n)$ hold with the following formulas 
      $$
      \begin{array}{ll}
      \lambda_{4n+1}=(-1)^{n+1}r_n^2, & \qquad \lambda_{4n+2}=(-1)^{n+1}(r_ns_n)^{-1},\\
      \lambda_{4n+3}=(-1)^{n+1}s_n^2,& \qquad \lambda_{4n+4}=(-1)^{n+1}(r_{n+1}s_n)^{-1}.
      \end{array} 
      $$
   From $(III)$, the link between the two sequences $(r_n)_{n\geq 1}$ and  $(s_n)_{n\geq 1}$ follows directly. Indeed, for $n\geq 1$, we have $s_n=r_n+r_{n+1}$. Furthermore, combining $(I)$ and $(III)$, with the values for $(\mu_4,\mu_8)$, we have
   $$r_{n+1}=2r_n+r_{n-1}\quad \text{for}\quad n\geq 2, \quad \text{with}\quad (r_1,r_2)=(12/25,32/25). \eqno{(R)}$$
   At last, we need to point at the link between the sequences $(r_n)_{n\geq 1}$ and $(\ell_n)_{n\geq 0}$.
   \newline For $n\geq 1$, we set $L_n=2\ell_n-\ell_{n-1}+1$. Then, it is easy to check that the sequence $(L_n)_{n\geq 1}$ satisfies the same second-order linear recursive relation $(R)$, with the initial conditions $(L_1,L_2)=(3,8)$. Therefore, we have the equality between the sequences $((25/4)r_n)_{n\geq 1}$ and $(L_n)_{n\geq 1}$.
      \newline So the proof of the theorem is complete.
   
 \vskip 1 cm
  
 \par We cannot conclude this note without asking the following question : are there other infinite words in $\Q$ defined in such a simple manner and leading to such a singular infinite continued fraction in $\Q((T^{-1}))$ ?

30 September 2024

\vskip 1 cm
\begin{thebibliography}{99}



\bibitem{1} B. Allombert, N. Brisebarre and A. Lasjaunias, \emph{On a two-valued sequence and related continued fractions in power series fields},
The Ramanujan Journal 45 (2018),  859--871.

\bibitem{2} W. Mills and D. P. Robbins, \emph{Continued fractions for certain algebraic power series.} J. Number Theory \textbf{23} (1986),
388--404.


\bibitem{3} A. Lasjaunias,\emph{ A survey of Diophantine approximation in fields of power series}, 
Monatshefte f\"{u}r Mathematik {\bf 130}  (2000), 211--229.


\bibitem{4} A. Lasjaunias and J-Y. Yao,\emph{ Hyperquadratic continued fractions and automatic sequences}, 
Finite Fields and their Applications {\bf 40}  (2016), 46--60.

\bibitem{5} A. Lasjaunias, \emph{Continued Fractions },   arXiv:1711.11276 (Novembre 2017). 



\bigskip
\end{thebibliography}
\end{document}